# SOME NOTES ON A SEMI-MARKOV MATRIX OCCURRING IN THE $M|M|\infty$ QUEUE PARAMETERS STUDY

Prof. Dr. **MANUEL ALBERTO M. FERREIRA**

Instituto Universitário de Lisboa (ISCTE – IUL), BRU - IUL, Lisboa, Portugal

manuel.ferreira@iscte.pt

**ABSTRACT**

The main objective of this work is to present a process to compute the Markov renewal matrix $R(t)$ for Markov renewal processes with countable infinite spaces, which semi-Markov matrixes $Q(t)$ are immigration and death type and assume a tridiagonal form. These processes occur often in practical applications. And the difficulty in obtaining friendly results for $R(t)$ is a great obstacle to its application in practical cases modelling. It is considered the application to the $M|M|\infty$ queue particular case.

## 1. INTRODUCTION

It is presented in the next section a process to compute the Markov renewal matrix $R(t)$ for Markov renewal processes, see (2), with countable infinite spaces, occurring often in practical applications, which semi-Markov matrixes $Q(t)$ are immigration and death type and assume a tridiagonal form.

Then, in the following section, this computation process is applied to a particular Markov renewal process that rules the $M|M|\infty$ queue, being its parameters designed in accordance.

The great motivation for this work is to get stochastic processes that approximate quite well infinite servers queues in order to obtain the exact values, or at least good approximations, for some interesting parameters difficult to obtain through other study methodologies, see for instance (3,4) and (7,8,9).

Despite the interesting and even fascinating calculations methods, the results obtained are quite disappointing. Actually they are only for the $R(t)$ entries Laplace transforms, assuming a very complicate form very hard, even impossible in many cases, to invert.

## 2. SEMI- MARKOVIAN TRIDIOGONAL TYPE MATRIX

In Markov renewal processes practical applications, occur often Markov renewal processes with a countable infinite states space which semi-Markov matrix is immigration and death type. That is

$$Q(t) = \begin{bmatrix} 0 & \tau_0(t) & 0 & & \\ \sigma_1(t) & 0 & \tau_1(t) & \cdots & \cdots \\ 0 & \sigma_2(t) & 0 & & \\ & \vdots & & \ddots & \vdots \\ & \cdots & & \cdots & \ddots \end{bmatrix}, t \geq 0 \quad (2.1)$$

where $\sigma_i(t)$ and $\tau_i(t)$ are any known functions of $t$ and $i$, and from the system parameters, that fulfill the following conditions:

- $0 \leq \sigma_i(t) \leq 1, i = 1,2,\ldots$
- $0 \leq \tau_i(t) \leq 1, i = 1,2,\ldots$
- $\sigma_i(t) + \tau_i(t) \leq 1, i = 1,2,\ldots$

and make the process regular.

The Markov renewal matrix $R(t)$ is defined as

$$R(t) = \sum_{n=0}^{\infty} Q^n(t) \quad (2.2)$$

and it is computed through its Laplace-Stieltjes transform, $\bar{R}(s)$, given by

$$\bar{R}(s) = \sum_{n=0}^{\infty} \bar{Q}^n(s) \quad (2.3)$$

where $\bar{Q}(s)$ is a matrix which entries are the $Q(t)$ entries Laplace-transforms. As $\sum_{n=0}^{\infty} \bar{Q}^n(s) = (I - \bar{Q}(s))^{-1}$ since $|I - \bar{Q}(s)| \neq 0$, the expression (2.3), and so (2.2), is equivalent to

$$\bar{R}(s)(I - \bar{Q}(s)) = I \quad (2.4)$$

being $I$ the identity matrix. As the Markov renewal process is a regular one, (2.4) has only one solution and it is equivalent to the equations system:

$$\sum_{k=0}^{\infty} [1 - \bar{\sigma}_k(s) - \bar{\tau}_k(s)] \bar{r}_{ik}(s) = 1, i = 0,1,\ldots; s \geq 0 \quad (2.5)$$

where it is settled that $\bar{\sigma}_0(0) = 0$. And, for each $i$, there is an infinite equations system

$$-\bar{\tau}_{j-1}(s)\bar{r}_{i,j-1}(s) + \bar{r}_{ij}(s) - \bar{\sigma}_{j+1}(s)\bar{r}_{i,j+1} = \delta_{ij}, j = 0,1,\ldots; s > 0 \quad (2.6)$$

where $\delta$ is the Kronecker delta and it is settled that $\bar{\tau}_{-1}(s) = 0$.

That is: to obtain the $\bar{r}_{ij}(s)$ it is necessary to solve the second order difference equations system (2.6) subject to the conditions (2.5), see (4). The possibility of obtaining explicit exact solutions depends on the $\bar{\tau}_j(s)$ and $\bar{\sigma}_j(s)$ functional forms. For instance, if

$$\bar{\tau}_j(s) = \frac{\rho}{j + \rho + \alpha s} \qquad (2.7)$$

and

$$\bar{\sigma}_j(s) = \frac{j}{j + \rho + \alpha s} \qquad (2.8)$$

with

$$\rho = \lambda \alpha \qquad (2.9)$$

that correspond to the Markov renewal processes that rules the $M|M|\infty$ queue, through the variable transformation, for a given $i$

$$\bar{t}_{ij}(s) = \frac{\alpha \bar{r}_{ij}(s)}{j + \rho + \alpha s}, j = 0,1,\dots \qquad (2.10)$$

and using generator functions it is possible to solve the infinite equations system (2.6), turning its solution in a non-constant coefficients first order differential equation, which solution is given in terms of the hypergeometric confluent series, as it will be seen in the next section.

### 3. THE $M|M|\infty$ QUEUE

Consider the matrix $Q(t)$ with $\tau_j(t)$ and $\sigma_j(t)$ such that $\bar{\tau}_j(s)$ and $\bar{\sigma}_j(s)$ are given by (2.7) and (2.8), respectively. Then (2.6) assumes the form

$$-\frac{\rho}{k-1+\rho+\alpha s}\bar{r}_{i,k-1}(s) + \bar{r}_{ik}(s) - \frac{k+1}{k+1+\rho+\alpha s}\bar{r}_{i,k+1}(s) = \delta_{ik},$$

$$k = 0,1,\dots; \bar{r}_{i,-1}(s) = 0, s \geq 0 \qquad (3.1).$$

Performing the variable change defined in (2.10), (2.1) becomes

$$-\rho \bar{t}_{i,k-1}(s) + (k + \rho + \alpha s)\bar{t}_{ik}(s) - (k+1)\bar{t}_{i,k+1}(s) = \delta_{ik}, k = 0,1,\dots (3.2).$$

Defining the generating functions $y_i(x)$, convergent for $|x| < 1$, as

$$y_i(x) = G_i(x;s) = \sum_{k=0}^{\infty} \frac{\alpha \bar{r}_{ik}(s)}{k+\rho+\alpha s} x^k = \sum_{k=0}^{\infty} \bar{t}_{ik}(s) x^k \qquad (3.3),$$

the expression (2.5) assumes the form

$$y_i(1) = G_i(1;s) = \frac{1}{s} \quad (3.4)$$

and multiplying in (3.2) the equation $k$ by $x^{k+1}, k = 0,1,\ldots$, summing all together member by member and using the generating functions defined in (3.3) the system resolution is converted in the resolution of the first order differential equation

$$(1-x)\frac{dy_i}{dx} - [\rho(1-x) + \alpha s]\, y_i(x) + \alpha x^i = 0, i = 0,1,\ldots \quad (3.5),$$

being (3.4) the initial condition.

Using the integrating factor

$$IF(x) = e^{-\rho x}(1-x)^{\alpha s - 1} \quad (3.6),$$

see for instance (5), equations (3.5) become exact differential equations. Then the solution is

$$y_i(x) = \alpha e^{-\rho(1-x)} \sum_{k=0}^{i} (-1)^k \binom{i}{k} \frac{(1-x)^k}{\alpha s + k} \phi\left(\alpha s + k, \frac{\alpha s + 1}{\alpha} + k + 1; \rho(1-x)\right),$$

$$i = 0,1,\ldots \quad (3.7),$$

obtained in terms of the hypergeometric confluent series: $\phi(\cdot,\cdot;\cdot)$.

Through formula (3.3),

$$\frac{\alpha n!}{n + \rho + \alpha s} \bar{r}_{in}(s) = \left(\frac{d^n y_i(x)}{dx^n}\right)_{x=0}$$

$$= \sum_{k=0}^{i} \binom{i}{k} \frac{\alpha^{1-k}}{\alpha s + k} \left(\frac{d^n}{dw^n}\left(w^k \phi(1, s\alpha + k + 1; w)\right)\right)_{w=-\rho} \left(\frac{dw}{dx}\right)^n_{x=0}, i$$

$$= 0,1,\ldots; n = 1,2,\ldots \quad (3.8)$$

where $w = -\rho(1-x)$. Finally, using Leibnitz's formula and the results for hypergeometric confluent series derivation, see (1),

$$\bar{r}_{in}(s) = \left(\frac{n}{\alpha} + \lambda + s\right) \sum_{k=0}^{i}(-1)^k \binom{i}{k} \frac{\alpha}{\alpha s + k} \sum_{j=0}^{\min(s,k)}(-1)^j \binom{k}{j}\left(\frac{\rho}{s\alpha + k + 1}\right)^{n-j} \phi(n$$

$$- j + 1, s\alpha + k + n - j + 1; -\rho), i, n = 0,1,\ldots; s \geq 0 \quad (3.9).$$

## CONCLUSIONS

Along this work, was evidenced that the Laplace transforms of the matrix $R(t)$ entries, $\bar{r}(s)$, are connected in a system of second order linear difference equations. For the particular case of exponentially, it was possible to find an explicit analytic solution to that system in terms of the hypergeometric confluent series. And there is no notice of any other situation for which this happens.